\documentclass[graybox]{svmult}

\usepackage{type1cm}        
\usepackage{makeidx}         
\usepackage{graphicx}        
\usepackage{multicol}        
\usepackage[bottom]{footmisc}

\usepackage{newtxtext}       %
\usepackage[varvw]{newtxmath}       
\usepackage{tikz-cd} 
\usepackage{tikz}
\usepackage{graphicx}
\usetikzlibrary{cd} 
\usetikzlibrary{calc}
\usetikzlibrary{decorations}
\usetikzlibrary{positioning}
\usetikzlibrary{shapes}
\usetikzlibrary{external}

\definecolor{symbols}{rgb}{0,0.7,1}
\definecolor{pagebackground}{rgb}{1,1,1}
\colorlet{connection}{red!80!black}
\colorlet{boxcolor}{blue!50}
\makeindex             
\makeatletter
\pgfdeclareshape{crosscircle}
{
	\inheritsavedanchors[from=circle] 
	\inheritanchorborder[from=circle]
	\inheritanchor[from=circle]{north}
	\inheritanchor[from=circle]{north west}
	\inheritanchor[from=circle]{north east}
	\inheritanchor[from=circle]{center}
	\inheritanchor[from=circle]{west}
	\inheritanchor[from=circle]{east}
	\inheritanchor[from=circle]{mid}
	\inheritanchor[from=circle]{mid west}
	\inheritanchor[from=circle]{mid east}
	\inheritanchor[from=circle]{base}
	\inheritanchor[from=circle]{base west}
	\inheritanchor[from=circle]{base east}
	\inheritanchor[from=circle]{south}
	\inheritanchor[from=circle]{south west}
	\inheritanchor[from=circle]{south east}
	\inheritbackgroundpath[from=circle]
	\foregroundpath{
		\centerpoint%
		\pgf@xc=\pgf@x%
		\pgf@yc=\pgf@y%
		\pgfutil@tempdima=\radius%
		\pgfmathsetlength{\pgf@xb}{\pgfkeysvalueof{/pgf/outer xsep}}%
		\pgfmathsetlength{\pgf@yb}{\pgfkeysvalueof{/pgf/outer ysep}}%
		\ifdim\pgf@xb<\pgf@yb%
		\advance\pgfutil@tempdima by-\pgf@yb%
		\else%
		\advance\pgfutil@tempdima by-\pgf@xb%
		\fi%
		\pgfpathmoveto{\pgfpointadd{\pgfqpoint{\pgf@xc}{\pgf@yc}}{\pgfqpoint{-0.707107\pgfutil@tempdima}{-0.707107\pgfutil@tempdima}}}
		\pgfpathlineto{\pgfpointadd{\pgfqpoint{\pgf@xc}{\pgf@yc}}{\pgfqpoint{0.707107\pgfutil@tempdima}{0.707107\pgfutil@tempdima}}}
		\pgfpathmoveto{\pgfpointadd{\pgfqpoint{\pgf@xc}{\pgf@yc}}{\pgfqpoint{-0.707107\pgfutil@tempdima}{0.707107\pgfutil@tempdima}}}
		\pgfpathlineto{\pgfpointadd{\pgfqpoint{\pgf@xc}{\pgf@yc}}{\pgfqpoint{0.707107\pgfutil@tempdima}{-0.707107\pgfutil@tempdima}}}
	}
}
\makeatother
\def\one{\mathbf{1}}

\tikzstyle{tinydots}=[dash pattern=on \pgflinewidth off \pgflinewidth]

\def\decorate#1#2{
	\ifnum#2>0
	\foreach \count in {1,...,#2}{
		let
		\p1 = (sourcenode.center),
		\p2 = (sourcenode.east),
		\n1 = {\x2-\x1},
		\n2 = {1mm},
		\n3 = {(1.3+0.6*(\count-1))*\n1},
		\n4 = {0.7*\n1}
		in 
		node[rectangle,fill=symbols,rotate=30,inner sep=0pt,minimum width=0.2*\n2,minimum height=\n2] at ($(sourcenode.center) + (\n3,\n4)$) {}
	}
	\fi
	\ifnum#1>0
	\foreach \count in {1,...,#1}{
		let
		\p1 = (sourcenode.center),
		\p2 = (sourcenode.east),
		\n1 = {\x2-\x1},
		\n2 = {1mm},
		\n3 = {(1.3+0.6*(\count-1))*\n1},
		\n4 = {0.7*\n1}
		in 
		node[rectangle,fill=symbols,rotate=-30,inner sep=0pt,minimum width=0.2*\n2,minimum height=\n2] at ($(sourcenode.center) + (-\n3,\n4)$) {}
	}
	\fi
}

\tikzset{
	dectriangle/.style 2 args={
		triangle,
		alias=sourcenode,
		append after command={\decorate{#1}{#2}}
	},
	dectriangle/.default={0}{0},
}

\tikzset{
	cross/.style={path picture={ 
			\draw[symbols]
			(path picture bounding box.south east) -- (path picture bounding box.north west) (path picture bounding box.south west) -- (path picture bounding box.north east);
	}},
	root/.style={circle,fill=green!50!black,inner sep=0pt, minimum size=1.2mm},
	dot/.style={circle,fill=pageforeground,inner sep=0pt, minimum size=1mm},
	dotred/.style={circle,fill=pageforeground!50!pagebackground,inner sep=0pt, minimum size=2mm},
	var/.style={circle,fill=pageforeground!10!pagebackground,draw=pageforeground,inner sep=0pt, minimum size=3mm},
	var2/.style={circle,fill=darkgreen,draw=pageforeground,inner sep=0pt, minimum size=3mm},
	kernel/.style={semithick,draw=green,shorten >=2pt,shorten <=2pt},
	kernels/.style={snake=zigzag,shorten >=2pt,shorten <=2pt,segment amplitude=1pt,segment length=4pt,line before snake=2pt,line after snake=5pt,},
	rho/.style={densely dashed,semithick,shorten >=2pt,shorten <=2pt},
	testfcn/.style={dotted,semithick,shorten >=2pt,shorten <=2pt},
	renorm/.style={shape=circle,fill=pagebackground,inner sep=1pt},
	labl/.style={shape=rectangle,fill=pagebackground,inner sep=1pt},
	xic/.style={very thin,circle,draw=symbols,fill=symbols,inner sep=0pt,minimum size=1.2mm},
	g/.style={very thin,rectangle,draw=symbols,fill=symbols!10!pagebackground,inner sep=0pt,minimum width=2.5mm,minimum height=1.2mm},
	xi/.style={very thin,circle,draw=symbols,fill=symbols!10!pagebackground,inner sep=0pt,minimum size=1.2mm},
	xies/.style={very thin,rectangle,fill=green!50!black!25,draw=symbols,inner sep=0pt,minimum size=1.1mm},
	xiesf/.style={very thin,rectangle,fill=green!50!black,draw=symbols,inner sep=0pt,minimum size=1.1mm},
	xix/.style={very thin,crosscircle,fill=symbols!10!pagebackground,draw=symbols,inner sep=0pt,minimum size=1.2mm},
	X/.style={very thin,cross,rectangle,fill=pagebackground,draw=symbols,inner sep=0pt,minimum size=1.2mm},
	xib/.style={thin,circle,fill=symbols!10!pagebackground,draw=symbols,inner sep=0pt,minimum size=1.6mm},
	xie/.style={thin,circle,fill=green!50!black,draw=symbols,inner sep=0pt,minimum size=1.6mm},
	xid/.style={thin,circle,fill=symbols,draw=symbols,inner sep=0pt,minimum size=1.6mm},
	xibx/.style={thin,crosscircle,fill=symbols!10!pagebackground,draw=symbols,inner sep=0pt,minimum size=1.6mm},
	kernels2/.style={very thick,draw=connection,segment length=12pt},
	keps/.style={thin,draw=symbols,->},
	kepspr/.style={thick,draw=connection,->},
	krho/.style={thin,draw=symbols,superdense,->},
	krhopr/.style={thick,draw=connection,superdense},
	triangle/.style = { regular polygon, regular polygon sides=3},
	not/.style={thin,circle,draw=connection,fill=connection,inner sep=0pt,minimum size=0.5mm},
	diff/.style = {very thin,draw=symbols,triangle,fill=red!50!black,inner sep=0pt,minimum size=1.6mm},
	diff1/.style = {very thin,dectriangle={1}{0},fill=red!50!black,draw=symbols,inner sep=0pt,minimum size=1.6mm},
	diff2/.style = {very thin,dectriangle={1}{1},fill=red!50!black,draw=symbols,inner sep=0pt,minimum size=1.6mm},
	diffmini/.style = {very thin,rectangle,fill=black,draw=black,inner sep=0pt,minimum size=0.75mm},
	kernelsmod/.style={very thick,draw=connection,segment length=12pt},
	rec/.style = {very thin,rectangle,fill=black,draw=black,inner sep=0pt,minimum size=2mm},
	cerc/.style={very thin,circle,draw=black,fill=symbols,inner sep=0pt,minimum size=2mm},
	stars/.style={very thin,star,star points=6,star point ratio=0.5, draw=black,fill=red,inner sep=0pt,minimum size=0.7mm},
	>=stealth,
}
\tikzset{
	root/.style={circle,fill=black!50,inner sep=0pt, minimum size=3mm},
	circ/.style={circle,fill=white,draw=black,very thin,inner sep=.5pt, minimum size=1.2mm},
	round1/.style={fill=white,outer sep = 0,inner sep=2pt,rounded corners=1mm,draw,text=black,thin,minimum size=1.2mm},
	circ1/.style={circle,fill=red!10,draw=red,very thin,inner sep=.5pt, minimum size=1.2mm},
	rect/.style={fill=white,outer sep = 0,inner sep=2pt,rectangle,draw,text=black,thin,minimum size=1.2mm},
	rect1/.style={fill=white,outer sep = 0,inner sep=2pt,rectangle,draw,text=black,thin,minimum size=1.2mm},
	round2/.style={fill=red!10,outer sep = 0,inner sep=2pt,rounded corners=1mm,draw,text=black,thin,minimum size=1.2mm},
	round3/.style={fill=blue!10,outer sep = 0,inner sep=2pt,rounded corners=1mm,draw,text=black,thin,minimum size=1.2mm}, 
	rect2/.style={fill=black!10,outer sep = 0,inner sep=2pt,rectangle,draw,text=black,thin,minimum size=1.2mm},
	dot/.style={circle,fill=black,inner sep=0pt, minimum size=1.2mm},
	dotred/.style={circle,fill=black!50,inner sep=0pt, minimum size=2mm},
	var/.style={circle,fill=black!10,draw=black,inner sep=0pt, minimum size=3mm},
	kernel/.style={semithick,draw=darkgreen},
	diag/.style={thin,shorten >=4pt,shorten <=4pt},
	kernel1/.style={thick},
	kernels/.style={snake=zigzag,shorten >=2pt,shorten <=2pt,segment amplitude=1pt,segment length=4pt,line before snake=2pt,line after snake=5pt},
	kernels1/.style={snake=zigzag,segment amplitude=0.5pt,segment length=2pt},
	rho1/.style={densely dotted,semithick},
	rho/.style={densely dashed,semithick,shorten >=2pt,shorten <=2pt},
	testfcn/.style={dotted,semithick,shorten >=2pt,shorten <=2pt},
	visible/.style={draw, circle, fill, inner sep=0.25ex},
	renorm/.style={shape=circle,fill=white,inner sep=1pt},
	labl/.style={shape=rectangle,fill=white,inner sep=1pt},
	xic/.style={very thin,circle,fill=symbols,draw=black,inner sep=0pt,minimum size=1.2mm},
	xi/.style={very thin,circle,fill=blue!10,draw=black,inner sep=0pt,minimum size=1.2mm},
	xib/.style={very thin,circle,fill=blue!10,draw=black,inner sep=0pt,minimum size=1.6mm},
	xie/.style={very thin,circle,fill=green!50!black,draw=black,inner sep=0pt,minimum size=1mm},
	xid/.style={very thin,circle,fill=symbols,draw=black,inner sep=0pt,minimum size=1.6mm},
	edgetype/.style={very thin,circle,draw=black,inner sep=0pt,minimum size=5mm},
	nodetype/.style={very thick,circle,draw=black,inner sep=0pt,minimum size=5mm},
	kernels2/.style={very thick,draw=connection,segment length=12pt},
	clean/.style={thin,circle,fill=black,inner sep=0pt,minimum size=1mm},	not/.style={thin,circle,fill=symbols,draw=connection,fill=connection,inner sep=0pt,minimum size=0.8mm},
	>=stealth,
}

\begin{document}

\title*{Resonances and computations}
\author{Yvain Bruned \and Frédéric Rousset \and Katharina Schratz }
\institute{Yvain Bruned \at University of Lorraine, Boulevard des Aiguillettes, 54506 Vandœuvre-lès-Nancy. \email{yvain.bruned@univ-lorraine.fr}	
	 \and Frédéric Rousset \at Université Paris- Saclay, 91405 Orsay Cedex, \email{frederic.rousset@universite-paris-saclay.fr}	
 \and Katharina Schratz \at Sorbonne University, 5 place Jussieu, 7505 Paris, \email{katharina.schratz@sorbonne-universite.fr}}
%
%
\maketitle
%

\abstract{The computation of  time dynamics arising in nonlinear time-dependent partial differential equations  is an ongoing challenge  in numerical analysis, especially once roughness comes into play. Classical numerical schemes in general fail to resolve the oscillatory behaviour in the solution which leads to numerical instabilities and loss of convergence. Dispersive equations, e.g., nonlinear Schr\"odinger, Korteweg--de Vries and wave equations, thereby pose in particular a big problem as in contrast to the parabolic setting, no strong smoothing can be expected, i.e., if the initial data is rough, the solution stays rough which makes their approximation a delicate task. In this review we give an overview on a new numerical ansatz which aims to tackle the time dynamics of nonlinear dispersive partial differential equations  even for very rough data. This is achieved by a resonance analysis and decorated tree formalism that draws its inpiration from the combinatorics used in the theory of regularity structures for solving singular SPDEs. One can hope to see this formalism applied in other contexts for dispersive PDEs and beyond.}

\section{Introduction}\label{sec:1}
The numerical computation of nonlinear partial differential equations (PDEs) for smooth solutions is nowadays in large parts well understood. A large toolbox of techniques to approximate effectively their time dynamics was established, involving splitting methods (with the intention of splitting the full equation into a series of simpler subproblems), exponential integrators (based on Duhamel's formula) and many more approaches. For an extensive overview of classical methods in context of dispersive PDEs we refer to \cite{EFHI09,Faou12,H2Tri,H2Tri,HochOst10,SanBook} and the references therein.

While classical techniques in many situations allow an efficient and well-suited numerical approximation, they in general fail to reproduce the behaviour of solutions once roughness comes into play. This is due to the fact that in general classical integration methods require a lot of regularity on the solutions in order to converge. If the solution is not regular enough, well-known theoretical error estimates break down and almost nothing is known on how to overcome this. Break down of convergence for rough data is not only a theoretical artefact, but also severely observed in numerical simulations. If the solution does not meet the required smoothness assumptions, the numerical solution becomes unstable and eventually blows up. More details on the struggle of classical methods for rough data are given in Section \ref{sec:2}.

The main question we want to address in this review is the following: How can we actually reproduce rough time dynamics in nonlinear PDEs an overcome classical smoothness assumptions?

To answer this question we will present a novel integration technique, so-called \emph{resonances as a computational tool}, with the aim to reduce as far as possible the high regularity assumptions needed in classical schemes. The novel idea of this  approach lies in building numerical schemes on the underlying structure of resonances of the PDE, instead of discretising the PDE itself. This is achieved by looking at iterations of Duhamel's formula and filtering out the dominant parts in the oscillations. To control the higher order oscillations, we use a decorated tree formalism with the frequencies as decorations. This will allow us to construct resonance based schemes up to arbitrary order.  Details are given in Section \ref{subsec:rc1} where we present the main ideas behind the derivation of low regularity schemes for dispersive equations. We also present some extensions that allow us to go deeper, e.g., preserving essential symmetries in the equation.

\subsection*{Acknowledgements}

{\small
	Y. B. gratefully acknowledges funding support from the European Research Council (ERC) through the ERC Starting Grant Low Regularity Dynamics via Decorated Trees (LoRDeT), grant agreement No.\ 101075208.
K.S. and F.R. gratefully acknowledges funding from the European Research Council (ERC) under the European Union’s Horizon 2020 research and innovation programme (grant agreement No. 850941).
}

\section{Classical methods and their problem for rough data}\label{sec:2}
To present the main challenges of classical methods for rough data, let us consider as a model problem the periodic Korteweg--de Vries equation
\begin{equation}\label{kdv}
\partial_t u(t,x) + \partial_x^3 u(t,x) = \frac12 \partial_x u^2(t,x) \qquad (t,x)\in \mathbb{R}\times \mathbb{T}
\end{equation}
with rough initial data
\begin{equation}\label{iniSig}
u(0,x) = u_0(x) \in H^\sigma(\mathbb{T}),
\end{equation}
where we want to choose $\sigma >0$ as small as possible (we will see later  how far we can actually push down the Sobolev index $\sigma$). Here, $H^\sigma(\mathbb{T})$ denotes the classical Sobolev space on   torus $\mathbb{T}$ with regularity $\sigma\geq 0$.

Before we present the main idea of {\em resonances as a computational tool}  -- a new method which ideally allows us to choose $\sigma$ in \eqref{iniSig} arbitrarily close to zero -- let us first give a brief overview of classical numerical methods  to solve the nonlinear PDE \eqref{kdv}. This will allow us to better understand the regularity assumptions of classical approaches  and get essential insight on how to overcome this.

\subsection{Splitting methods} One of the most famous  numerical methods to approximate the time dynamics of KdV equation \eqref{kdv}  and in general nonlinear dispersive models are splitting methods (see, e.g., \cite{Faou12,H2Tri,HLRS10,HUV,JL00,Lubich08,McLacQ02,Shi}). The main idea thereby lies in splitting the full equation into a series of simpler subproblems. One then solves these subproblems (either exactly or with a numerical scheme) and composes the sub-flows to obtain an approximation to the flow of the original -- full -- equation. 

 Looking at the KdV equation \eqref{kdv} we face two main challenges numerically: The differential operator $\partial_x^3$ and the non-linear Burger's term  $\frac12 \partial_x u^2$. The idea of splitting is now the following: Instead of solving the full problem \eqref{kdv} we consider its linear and nonlinear part  separetely, i.e., 
\begin{align}
\text{ (L) } \quad \partial_t u(t,x) = - \partial_x^3 u(t,x) \quad \text{and} \quad\text{ (N) } \quad \partial_t u(t,x) = \frac12 \partial_x u^2(t,x).
\end{align}
The  main advantage thereby lies in the fact that the time dynamics in both subproblems can be solved exactly: (L) in Fourier space (in case of a spatial discretisation we can couple this with for instance a Fourier pseudo spectral method) and (N) with the method of characteristics. 

\begin{remark}
In the end we are of course interested in a full discretisation as we have a PDE both in time and in space. This  will introduce an additional spatial error. In case of splitting for KdV we will have a  spectral error for the linear problem (L) as well as an interpolation error for the nonlinear problem (N)  as we do not necessarily end up on the grid points any more when using the method of characteristics.  { The latter requires an interpolation, e.g., Lagrange interpolation, or the use of WENO schemes, e.g., upwind scheme, see for instance \cite{Lagi,WENO}.}
 In this review we only focus on the time dynamics and hence neglect the additional error introduced by the spatial discretisation of PDEs. For fully discrete resonance based schemes in context of Navier--Stokes equation and Schr\"odinger equations we refer for instance to \cite{LiS,ROS24}.
\end{remark}
A simple composition of the kinetic flow $\phi_L^t$ and nonlinear flow  $\phi_N^t$  leads at first-order  to the so-called first order Lie splitting method (with time step size $\tau$)
\begin{equation}\label{Lie}
\phi_L^\tau \circ \phi_N^\tau
\end{equation}
and its symetrised version the classical second-order Strang splitting
\begin{equation}\label{StrangKdV}
\phi_L^{\tau/2} \circ \phi_N^\tau\circ \phi_L^{\tau/2} .
\end{equation}
We refer to an overview of splitting methods to \cite{H2Tri,HUV,McLacQ02,Shi}.

Here we are in   interested in the local error structure the splitting approach introduces for KdV, as this will allow us to answer the question what kind of regularity assumptions this imposes on the solution and more importantly, how small can we actually choose the regularity $\sigma$ in the initial value \eqref{iniSig} in case of solving KdV with Strang splitting?

First rigorous error estimates where achieved in  \cite{HKRT12}. { Estimates requiring less regularity on the solution have been established in \cite{HLR12} showing that} the local error of  the Strang splitting \eqref{StrangKdV} for the KdV equation \eqref{kdv} is driven by the triple commutator 
\begin{equation}\label{Tcomm}
[ \partial_x^3,[ \mathcal{N}, \partial_x^3]] 
\end{equation}
where $\mathcal{N}$ denotes the Burgers nonlinearity $\mathcal{N}(u) = \frac12 \partial_x u^2 $. In this commutator structure the  ``worst-order'' derivatives $\partial_x^7$ and $\partial_x^6$  drop  such that formally $[ \partial_x^3,[ \mathcal{N}, \partial_x^3]] $ is a fifth-order differential operator  (cf. \cite[Lemma 6.1]{HLR12}) which implies  a local error behaviour  at order
$$
\mathcal{O}\left(  \tau^3 \partial_x^5 u  \right).
$$
Hence, in order to reach second-order convergence Strang splitting requires at least the boundedness of five additional derivatives of the solution.


\subsection{Exponential integrators} 
Another well-known method to solve KdV-type equations are exponential integrators  (see, e.g., \cite{CCO08,Klein06,HochOst10,OstS19} and the references therein). The main idea  thereby lies in discretising   Duhamel's formula which for KdV \eqref{kdv}  takes the form
\begin{equation}\label{Duh}
u(t) =e^{- t \partial_x^3} u(0) +\frac12 \partial_x e^{- t \partial_x^3} \int_0^t e^{ s \partial_x^3}  u^2(s)\mathrm{d}s.
\end{equation}
At  time $t_{n+1}=t_n + \tau$ we find (considering Duhamels formula on the time interval $[0,\tau]$  with  initial value $u(t_n)$) that
 \begin{equation}\label{DuhN}
u(t_{n+1}) =e^{- t \partial_x^3} u(t_n) +\frac12 \partial_x e^{- \tau \partial_x^3} \int_0^\tau e^{ s \partial_x^3}  u^2(t_n+s)\mathrm{d}s.
\end{equation}
Exponential integrator schemes are based on  Taylor series expansion of the solution  within the integral. At first order they build on the first-order Taylor series expansion 
\begin{equation}\label{taylor1exp}
u(t_n+s) = u(t_n) + \mathcal{O}( s u').
\end{equation}
Plugging the approximation \eqref{taylor1exp} into Duhamel's formula \eqref{DuhN} leads to the first-order exponential integrator scheme
\begin{equation}\label{expInt}
u^{n+1} = e^{- \tau \partial_x^3} u^n + \frac12\tau \partial_x \varphi_1(-\tau \partial_x^3)   (u^n)^2 \quad \text{with}\quad\varphi_1(z) = \frac{e^{z}-1}{z}.
\end{equation}
From the Taylor series expansion \eqref{taylor1exp} we easily see that the local error of the first-order exponential integrator method \eqref{expInt} is driven by the time derivative $u'$, where by the KdV equation \eqref{kdv} in the sense of regularity we have 
\begin{equation}\label{errExp}
\mathcal{O}(u') = \mathcal{O}(\partial_x^3 u).
\end{equation}
Hence,  first-order convergence requires the boundedness of at least three additional derivatives. Higher-order exponential integrators are based on higher-order Taylor series expansion of the solution within the integral, i.e., 
at second-order one takes
\begin{equation}\label{taylor2exp}
u(t_n+s) = u(t_n) +s u'(t_n) + \mathcal{O}(s^2 u'')
\end{equation}
and replaces the time derivative $u'(t_n)$ by the equation itself $u'(t_n) = - \partial_x^3 u(t_n)+ \frac12 \partial_x u^2(t_n)$.   Due to the local error scaling
$$
 \mathcal{O}(u'')= \mathcal{O}(\partial_x^6 u)
$$
we see that second-ordre exponential integrator method requires the boundedness of at least six additional derivatives.
\subsection{Structure of the solution and classical methods}
From the  local error structure \eqref{Tcomm} and~\eqref{errExp} we see that both splitting and exponential integrator methods require sufficiently smooth solutions to guarantee the boundedness of their local error.  A natural question therefore arises:\\{\em What happens for rough data? And can we actually construct numerical schemes which allow  convergence for rougher data than classical methods do?}\smallskip

To answer this question we first have to understand the underlying structure of the solution to KdV \eqref{kdv}   in a better way. For this purpose we turn back to Duhamel's formula \eqref{Duh}. However, instead of approximating its solution $u(s)$ by a classical Taylor series expansion (as done in exponential integrator schemes) we look at iterations of   Duhamel's formula: Using that
$$
u(s) = e^{- s \partial_x^3}u(0) + \int_0^s \ldots \mathrm{d}s_1
$$
we find
\begin{equation}\label{DuhIt}
u(t) = e^{-t\partial_x^3} u(0) +\frac12 \partial_x e^{-t\partial_x^3} \int_0^t e^{ s\partial_x^3} \left(e^{ -s\partial_x^3}u(0)  \right)^2 \mathrm{d}s +  \int_0^t  \int_0^s \ldots \mathrm{d}s_1\mathrm{d}s.
\end{equation}
If we want to get a rough idea of the dynamics of the solution, we can at first forget about the higher order iterations, and neglect the double integral $ \int_0^t  \int_0^s\mathrm{d}s_1\mathrm{d}s$ in \eqref{DuhIt}. We then see that the underlying structure of the solution $u(t)$ is driven by the nonlinear frequency interaction of $\partial_x^3$ and $-\partial_x^3$ with leading oscillations
\begin{equation}\label{osc}
\mathrm{Osc}(s,\partial_x^3, u(0)) = \partial_x e^{ s\partial_x^3} \left(e^{ -s\partial_x^3}u(0)  \right)^2.
\end{equation}
Numerical schemes stay close to the structure of the solution    if they   resolve  the leading oscillations~\eqref{osc}. 

A closer look, however, shows that splitting methods and exponential integrators in general neglect the nonlinear frequency interactions in \eqref{osc}: Lie splitting \eqref{Lie}  is based on the frequency approximation
 \begin{equation}\label{Lie:freq}
\mathrm{Osc}(s,\partial_x^3, u(0)) \approx \partial_x u^2(0)
\end{equation}
while exponential integrator methods  swallow all frequencies within the nonlinearity based on  \begin{equation}\label{Exp:freq}
\mathrm{Osc}(s,\partial_x^3, u(0)) \approx \partial_x 
e^{ s\partial_x^3}u^2(0).
\end{equation}
In case of smooth solutions, for which $\partial_x^3 u$   is well defined in the space of interest   linearisation of frequencies such as \eqref{Lie:freq} and \eqref{Exp:freq}   in general lead to good approximations of the exact solution \eqref{DuhIt}. This can be seen by a simple Taylor series expansion of the oscillations
\begin{equation}\label{taylor:class}
e^{s\partial_x^3 } u = u + \mathcal{O}(s \partial_x^3 u).
\end{equation}
Expansion \eqref{taylor:class} introduces a small remainder of order $s$ as long as $u$ is sufficiently smooth, i.e., $ \partial_x^3 u$ is bounded.

 For rough solutions, on the other hand, for which 
$$
 \partial_x^3  u
$$
is unbounded,   approximations such as \eqref{Lie:freq} and \eqref{Exp:freq}   break down as the linearisation of the frequency interactions \eqref{taylor:class} is no longer valid. The latter is not only a theoretical artefact stemming from error analysis, but also drastically observed in numerical experiments, where we observe instability and loss of convergence for rough data.\smallskip

The main aim of {resonances as a computational tool} is to  overcome this by stepping away from linearising the frequency interactions as done by classical methods  (e.g., \eqref{Lie:freq} and \eqref{Exp:freq}) towards new schemes which deeply embed the nonlinear frequency interactions in PDEs  into the numerical discretisation. In general this allows us to achieve reliable approximations for much rougher data than classical schemes can handle.  

In Section \ref{subsec:rc1} we will explain this idea of tackling nonlinear frequency interactions in a nonlinear way in detail on the concrete example of  periodic KdV~\eqref{kdv}. 

We refer to \cite{HS16,OS18,SWZ20} for  first results on resonance based schemes for KdV, NLS and Dirac equations and the recent review article on first- and second-order methods \cite{RS24}.

\section{Resonances as a computational tool: The main idea}\label{subsec:rc1}

 We rewrite in Fourier space the Duhamel's formula given in \eqref{Duh} and we obtain the following equation for the $k$-th Fourier coefficient  of the solution:
\begin{equation} \label{duhamel_Fourier}
	\begin{aligned}
		u_k(t) & = e^{i t k^3 }   v_k  + i k \sum_{\substack{k_1,k_2 \in \mathbb{Z}\\k_1+k_2 = k} }  e^{ i t k^3}  \int_0^t e^{-is k^3} 	u_{k_1}(s) 	u_{k_2}(s) ds,
	\end{aligned}
\end{equation}
where the linear operators $\partial_x $, $ e^{- t \partial_x^3} $, $ e^{ t \partial_x^3} $ are sent respectively in Fourier space to $ik $, $ e^{i t k^3} $, $ e^{ - it k^3} $. The splitting of $ k $ into $ k_1, k_2$ comes from the fact that moving to Fourier space, pointwise products become convolution products on the frequency.
We fix $r=0$ and iterate \eqref{duhamel_Fourier} by replacing $ u_{k_j}(t) $ by
\begin{equation*}
	u_{k_j}(t) =  e^{i t k_j^2 }  v_{k_j} +  \mathcal{O}(t),
\end{equation*}
with $ j \in \lbrace 1,2 \rbrace  $.
We obtain the following first order approximation of the $k$-th Fourier coefficient $ u_k(t) $:
\begin{equation*}
	\begin{aligned}
	u_k(t) & = e^{ i t  k^3}  v_k   + i k \sum_{\substack{k_1,k_2 \in \mathbb{Z}\\k_1+k_2 = k} }  e^{ i t k^3}  \int_0^{t} e^{- is k^3} (e^{i s k_1^3} v_{k_1}) (e^{ i s  k_2^3}  v_{k_2})  ds  + \mathcal{O}(t^2).
	\end{aligned}
\end{equation*}
One can encode the previous Duhamel iterates using a decorated tree series. We denote by $ U_k^r(v,t) $ the first iterated integrals of size $r$ of the Duhamel expansion. These are integrals with $r$ integrations in time. One has
\begin{equation*}
	|u_k(t)  - U^r_k(v,t)| = \mathcal{O}(t^{r+1})
	\end{equation*}
where the regularity asked on the initial data hidden in the notation $\mathcal{O}$ corresponds to the regularity needed to 
define the first iterated integrals up to order $r$. Decorated trees are used to provide a precise description of $ U^r_k(v,t) $ and to implement in a systematic way the resonance analysis. One uses the following B-series type formula
\begin{equation}
	\label{scheme_equation}
	 U^r_k(v,t) = \sum_{T \in \mathcal{T}_k^r} \frac{\Upsilon(T)(v)}{S(T)} (\Pi T)(t)
\end{equation} 
where $\mathcal{T}_k^r$ is a suitable set of decorated trees of size $r$, $S(T)$ is a symmetry factor and $\Upsilon(T)$ is an elementray differential associated with $ T $ depending on the initial data $v$. The map $ \Pi $ sends $T$ to an oscillatory integral.
Such a formalism is reminiscence of the B-series { (named
	after Butcher)} that describes numerical schemes for ODEs and PDEs. We refer to \cite{Butcher72,MR2657947,H2Tri} { for references on B-series }   and to  \cite{Christ, oh1, Gub11,LO13,HLO20} for references on
	tree series used in the context of dispersive equations . The formalism \eqref{scheme_equation} is largely inspired from  the treatment of  singular stochastic partial differential equations (SPDEs) via regularity structures in \cite{reg,BHZ,BCCH} where a local ansatz for the solutions takes a similar form. It is also connected to rough paths \cite{Lyo98,Gub04,Gub10}, one of the foundational ideas of regularity structure. {Note that the 
	main difference with the classical B-series comes from the decorations 
	and the application from set of trees to differential objects. Indeed, the decorated trees encode iterated integrals $ (\Pi T)(t)$ but also elementary differentials $\Upsilon(T)(v)$.} 
Below, we describe the series for $r=2$. One has
\begin{equation*}
	\mathcal{T}^{2}_k  = \lbrace  T_0,  T_1, T_2, \, \, k_1, k_2, k_3 \in \mathbb{Z} \rbrace, 
\end{equation*}
where the decorated trees $ T_i $ are given by
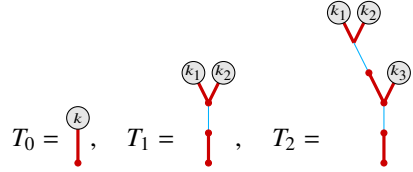
\begin{equation*}
	T_0  = \begin{tikzpicture}[scale=0.2,baseline=-5]
		\coordinate (root) at (0,1);
		\coordinate (tri) at (0,-2);
		\draw[kernels2] (tri) -- (root);
		\node[var] (rootnode) at (root) {\tiny{$ k $}};
		\node[not] (trinode) at (tri) {};
	\end{tikzpicture}, \quad  T_1   =   \begin{tikzpicture}[scale=0.2,baseline=-5]
		\coordinate (root) at (0,2);
		\coordinate (tri) at (0,0);
		\coordinate (tri1) at (0,-2);
		\coordinate (t1) at (-1,4);
		\coordinate (t2) at (1,4);
		\draw[kernels2] (t1) -- (root);
		\draw[kernels2] (t2) -- (root);
		\draw[kernels2] (tri) -- (tri1);
		\draw[symbols] (root) -- (tri);
		\node[not] (rootnode) at (root) {};
		\node[not] (trinode) at (tri) {};
		\node[not] (trinode) at (tri1) {};
		\node[var] (rootnode) at (t1) {\tiny{$ k_1 $}};
		\node[var] (trinode) at (t2) {\tiny{$ k_2 $}};
	\end{tikzpicture}, \quad T_2 = \begin{tikzpicture}[scale=0.2,baseline=-5]
		\coordinate (root) at (0,2);
		\coordinate (tri1) at (0,-2);
		\coordinate (tri) at (0,0);
		\coordinate (t1) at (-1,4);
		\coordinate (t11) at (-2,6);
		\coordinate (t12) at (-3,8);
		\coordinate (t13) at (-1,8);
		\coordinate (t2) at (1,4);
		\draw[kernels2] (t11) -- (t13);
		\draw[kernels2] (t11) -- (t12);
		\draw[kernels2] (t1) -- (root);
		\draw[kernels2] (tri) -- (tri1);
		\draw[symbols] (t1) -- (t11);
		\draw[kernels2] (t2) -- (root);
		\draw[symbols] (root) -- (tri);
		\node[not] (rootnode) at (root) {};
		\node[not] (trinode) at (tri) {};
		\node[not] (trinode) at (tri1) {};
		\node[not] (trinode) at (t1) {};
		\node[var] (rootnode) at (t12) {\tiny{$ k_{\tiny{1}} $}};
		\node[var] (rootnode) at (t13) {\tiny{$ k_{\tiny{2}} $}};
		\node[var] (trinode) at (t2) {\tiny{$ k_{\tiny{3}} $}};
	\end{tikzpicture}
\end{equation*} 
where for $T_1$, one has $ k  =k_1 + k_2 $ and for $T_2$ one has $k= k_1 + k_2 + k_3$. For the decorated trees written above, we have used the following coding:  An edge  $ \begin{tikzpicture}[scale=0.2,baseline=-5]
	\coordinate (root) at (0,-1);
	\coordinate (tri) at (0,0.5);
	\draw[kernels2] (root) -- (tri);
\end{tikzpicture} $ corresponds to a factor $ e^{i t k^3} $, while an edge  $  \begin{tikzpicture}[scale=0.2,baseline=-5]
\coordinate (root) at (0,-1);
\coordinate (tri) at (0,0.5);
\draw[symbols] (root) -- (tri);
\end{tikzpicture}  $  corresponds to an integral $ i k \int^{t}_0 e^{-i s k^3} \cdots d s $. The decorations on the leaves correspond to the frequencies that add up in the intern node up to $k$ which is the decoration of the node connected to the root of the tree. We have omitted these node decorations. One has
\begin{equation*}
	\begin{aligned}
	(\Pi T_0)(t)  & = e^{it k^3}, \quad 	(\Pi T_1)(t)  = i k e^{it k^3} \int^t_0  e^{-is (k^3-k_1^3 -k_2^3) }   ds \\ (\Pi T_2)(t)   & = - k k_{12} e^{it k^3} \int^t_0  e^{-is (k^3-k_{12}^3 -k_3^3) } v \int_0^s  e^{-i\bar{s} (k_{12}^3-k_{1}^3 -k_2^3) }  d \bar{s}  ds
	\end{aligned}
\end{equation*}
where $k_{12} = k_1 + k_2$. The symmetry factor $S(T)$ corresponds to the number of internal symmetries of the tree $T$ taking the edge decorations into account but not the node decorations. One obtains
\begin{equation*}
	S(T_0) = 1, \quad S(T_1) = 2, \quad S(T_2) = 2.
\end{equation*}
The $2$ in $ S(T_1) $ and $S(T_2)$ comes from the fact that one can permute the two leaves decorated by $k_1$ and $k_2$. The elementary differential $ \Upsilon(T)(v) $ corresponds to a product of initial data associated with the leaves of a decorated tree. One has to take into account also a factor connected to the structure of the tree: In the case of the KdV equation, this factor is $ 2^m $ where $m$ is the number of inner nodes in the tree. One gets
\begin{equation*}
	 \Upsilon(T_0)(v) = v_k, \quad \Upsilon(T_1)(v) = 2 v_{k_1} v_{k_2}, \quad \Upsilon(T_2)(v) = 4 v_{k_1} v_{k_2} v_{k_3}.
\end{equation*}
Let us remark that the sum over elements $ \mathcal{T}_k^r $ is infinite. One has a finite number of tree shapes and the infinity comes from the frequencies $k_1, k_2, k_3$ on which one performs an infinite summation. The term $U^r_k(t)$ gives a first approximation of order $r+1$ but one has to discretise each $(\Pi T_i)(t)$ such that its discretisation is easy to rewrite in physical space via usual differential operators. To produce a better scheme, one has to embed the resonance into the discretisation. Let us consider the first iterated integral
\begin{equation*}
	(\Pi T_1)(t)  = i k e^{it k^3} \int^t_0  e^{-is (k^3-k_1^3 -k_2^3) }   ds, \quad k= k_1 + k_2.
\end{equation*}
One observes the following
\begin{equation}
\label{factorization}
	k^3-k_1^3 -k_2^3 = 3 k_1 k_2 (k_1 + k_2) = \mathcal{L}_{\text{\tiny{low}}}.
\end{equation}
We use the notation $ \mathcal{L}_{\text{\tiny{low}}} $ with the subscript $ \tiny{\text{low}} $  to stress that the term obtained is of order lower than the operator meaning that one has at most a factor  $ k_i^2 $ which is lower than $k_i^3$.
Then, by exact integration when $  k_i \neq 0$ and $k_1 + k_2 \neq 0$, one gets
\begin{equation*}
	\int^t_0  e^{-is  \mathcal{L}_{\text{\tiny{low}}}}   ds = 
	- \frac{e^{-it  \mathcal{L}_{\text{\tiny{low}}}}-1}{i 3 k_1 k_2 (k_1 + k_2)}.
	\end{equation*}
One can notice that the factor at the denominator in the previous term can be easily mapped back to physical space. { The same is true for $ e^{-it  \mathcal{L}_{\text{\tiny{low}}}} $}. Indeed, one has
\begin{equation*}
	\frac{1}{k_1 k_2 (k_1 + k_2)} v_{k_1} v_{k_2} \longrightarrow
	\partial_x^{-1} \left(  \partial_x^{-1} v \right)^2
\end{equation*}
{ and by \eqref{factorization} (using it from the right to the left)
$$
e^{-it  \mathcal{L}_{\text{\tiny{low}}}} v_{k_1} v_{k_2} = e^{-it  k^3} e^{it  k_1^3}e^{it  k_2^3} v_{k_1} v_{k_2}
 \longrightarrow 
 e^{t \partial_x^3} \left( e^{-t \partial_x^3} v \right)^2
.
$$}
\noindent This is a crucial observation as the first iterated integral can be computed exactly without any further regularity on the initial data. Therefore, if we denote by $(\Pi^r T_1)(t)$ the approximation of order $r$ of $(\Pi T_1)(t)$ then one has
\begin{equation*}
	(\Pi^r T_1)(t) =  (\Pi T_1)(t) = 	- i k_{12} e^{it k^3} \frac{e^{-it  \mathcal{L}_{\tiny{\text{low}}}}-1}{i 3 k_1 k_2 (k_1+k_2)}, 
\end{equation*}
Now, we have to face the discretisation of the second order integral $ (\Pi T_2)(t)$. One can maybe hope for the same type of trick but it is not the case. Indeed, if we substitute the previous approximation, one has
\begin{equation*}
	(\Pi T_2)(t)    =  k k_{12} e^{it k^3} \int^t_0  e^{-is (k^3-k_{12}^3 -k_3^3) }  \frac{e^{-is  \mathcal{L}_{\tiny{\text{low}}}}-1}{i 3 k_1 k_2 (k_1+k_2)}  ds.
\end{equation*}
Then, one has two contributions to the resonance analysis which are given by
\begin{equation}
	\label{two_cases}
	\begin{aligned}
	k^3-k_{12}^3 -k_3^3  & = 3 k_{12} k_3 k,
	\\ k^3-k_{12}^3 -k_3^3 + \mathcal{L}_{\text{\tiny{low}}} & = 
	 k^3-k_{1}^3 - k_2^3 -k_3^3.  
	 \end{aligned}
	\end{equation}
The first term in the computation above can lead to an exact integration. Unfortunatly, the second term does not have a nice factorisation. In this case, we have to Taylor expand everything which is still better than classical schemes. As we want a very simple scheme, we do not perform any exact integration and Taylor expand everything as the local error is the same.
This gives the following approximation:
\begin{equation*}
	\begin{aligned}
	(\Pi T_2)(t)   &=  - k k_{12} e^{it k^3} \int^t_0  e^{-is (k^3-k_{12}^3 -k_3^3) }  \int_0^s \left(   1 + \mathcal{O}(\bar{s} (k_1^2 k_2 + k_{2}^2k_1)) \right)  d \bar{s}  ds
\\ 	= &- k k_{12} e^{it k^3} \int^t_0 s \left( 1 +  \mathcal{O}(s (k_{12}^2 k_3 + k_{3}^2k_{12}))  \right)  ds +  \mathcal{O}(t^3 k k_{12} (k_1^2 k_2 + k_{2}^2k_1))
\\ 	= &- k k_{12} e^{it k^3} \frac{t^2}{2} +    \mathcal{O}(t^3  k k_{12} (k_{12}^2 k_3 + k_{3}^2k_{12}))  +  \mathcal{O}(t^3 k k_{12} (k_1^2 k_2 + k_{2}^2k_1))
	\end{aligned}
\end{equation*}
where the local error term with a polynomial in the frequency inside the $\mathcal{O}(\cdot)$ means that we have to have enough regularity on the initial data $v$ in some suitable Sobolev space. One can observe that the two error terms ask the same regularity on the initial data which is $5$ derivatives. We therefore set
\begin{equation*}
	(\Pi^2 T_2)(t) = - k k_{12} e^{it k^3} \frac{t^2}{2} 
\end{equation*}
and one has
\begin{equation*}
\left|	\frac{\Upsilon(T)(v)}{S(T)} (\Pi T_2)(t) - \frac{\Upsilon(T)(v)}{S(T)} (\Pi^2 T_2)(t) \right| = \mathcal{O}( t^3 P(k_1,k_2,k_3) \Upsilon(T)(v) )
\end{equation*}
where $ P $ is a polynomial with degree at most equal to $5$ in each variable $k_i$.
Computing by hand without any structure for higher-oder schemes is difficult. This is where the decorated trees' formalism is useful providing a framework that organises these computations. 

One important remark is that in general in the computation of the resonance, one does not expect to see a decrease in the order of the polynomial. Indeed,
in \eqref{two_cases}, we can observe that one moves degree threee to degree two. This is not the case for the Schrödinger equation and other dispersive equations. In general, one has the following structure:
\begin{equation*}
\mathcal{L} =	P(k)  + \sum_{i=1}^n a_i P(k_i)
\end{equation*} 
where the $a_i$ belongs to $ \lbrace -1, 1 \rbrace $ and the polynomial $P$ is the differential operator of the equation translated into Fourier space. One expects a decomposition of the form
\begin{equation}
	\label{good_splitting}
		P(k)  + \sum_{i=1}^n a_i P(k_i) =  \mathcal{L}_{\text{\tiny{dom}}} + \mathcal{L}_{\text{\tiny{low}}}
\end{equation}
where $\mathcal{L}_{\text{\tiny{dom}}}$ contains the higher monomials in some frequency $k_i$. If $ \mathcal{L}_{\text{\tiny{dom}}} = 0 $ then we are in a case similar to the KdV equation. Otherwise, we want $\mathcal{L}_{\text{\tiny{dom}}}$ to be in a specific form for being able to map $ 	\int_0^t e^{i s \mathcal{L}_{\text{\tiny{dom}}}} ds $ back into physical space with usual differential operators. This is the case for $ \mathcal{L}_{\text{\tiny{dom}}} = k_i^n $ or any linear combination of these monomials with coefficients in $ \lbrace -1, 1 \rbrace $. Then, one proceeds as follows
\begin{equation} \label{resonance_decomposition}
	\begin{aligned}
	\int_0^t e^{i s \mathcal{L}} ds &= 	\int_0^t e^{i s \mathcal{L}_{\text{\tiny{dom}}}} e^{i s \mathcal{L}_{\text{\tiny{low}}}}ds
	\\ &= 	\int_0^t e^{i s \mathcal{L}_{\text{\tiny{dom}}}} (1 + \mathcal{O}(s\mathcal{L}_{\text{\tiny{low}}})) ds = \int_0^t e^{i s \mathcal{L}_{\text{\tiny{dom}}}} ds +  \mathcal{O}(s^2\mathcal{L}_{\text{\tiny{low}}})
	\end{aligned}
	\end{equation}
which gives a discretisation requiring a regularity given by  $ \mathcal{L}_{\text{\tiny{low}}} $ instead of $ \mathcal{L} $. The approximation above could be seen as Taylor expanding the lower part of the operator $\mathcal{L}$ denoted by $\mathcal{L}_{\text{\tiny{low}}}$ and integrating exactly its dominant part $\mathcal{L}_{\text{\tiny{dom}}}$. Let us stress that in the deterministic setting, one has an exact expression for the integral  $ \int_0^t e^{i s \mathcal{L}_{\text{\tiny{dom}}}} ds $ 
In a stochastic context, one has to deal with stochastic iterated integrals as in \cite{AGB23}. For those, one does not have pathwise calculus for these iterated integrals and therefore, one is limited to a few iterations for a low regularity scheme.

The splitting observed in \eqref{good_splitting} happens in fact in the context of the nonlinear cubic Schrödinger equation (NLS) given for $ (t,x) \in \mathbb{R}_+ \times \mathbb{T}^d $ by
\begin{equation} \label{equa_example}
	i \partial_t u + \Delta  u = \vert u\vert^2 u, \quad u(0,x) = v(x).
\end{equation}
The Duhamel formula for the $k$-th Fourier coefficient of the solution is given by
\begin{equation*} 
	\begin{aligned}
		u_k(t) & = e^{- i t k^2 }   v_k -\sum_{\substack{k_1,k_2,k_3 \in \mathbb{Z}^d\\-k_1+k_2+k_3 = k} } i e^{ -i t k^2} \int_0^t e^{i s k^2} \overline{	u_{k_1}}(s) 	u_{k_2}(s)	u_{k_3}(s)ds.
	\end{aligned}
\end{equation*}
{ Now as we work for $d \geq  1$, the notations $k^2$ and $k \ell$ mean for $k, \ell \in \mathbb{Z}^d$
$$
k^2 = \sum_{i=1}^d (k^{(i)})^2, \quad k \ell = \sum_{i=1}^d k^{(i)} \ell^{(i)}
$$
where $k = (k^{(1)},...,k^{(d)})$ and $\ell = (\ell^{(1)},...,\ell^{(d)})$.}
The first iterated integral  is 
\begin{equation} \label{first_integral_NLS}
- i e^{ -i t k^2} \int_0^t e^{i s k^2} e^{i s k_1^2} e^{-i s k^2_2} e^{-i s k_3^2} 	ds \, \overline{v}_{k_1} v_{k_2} v_{k_3}
\end{equation}
where $ k=-k_1 + k_2 + k_3 $.  Then, the resonance inside the integral is given by
\begin{equation*}
	k^2 + k_1^2 - k_2^2 - k_3^2 =  2 k_1^2 - 2 k_1 (k_2 + k_3) + 2 k_2 k_3.
\end{equation*}
In this case, one has $ \mathcal{L}_{\text{\tiny{dom}}} \neq 0 $ and
\begin{equation*}
	\mathcal{L}_{\text{\tiny{dom}}}  = 2 k_1^2, \quad 	\mathcal{L}_{\text{\tiny{low}}}  = - 2 k_1 (k_2 + k_3) + 2 k_2 k_3.
\end{equation*}
We are therefore in the case where exact integration and full Taylor expansion are needed at the same time. The discretisation of \eqref{first_integral_NLS} is given by
\begin{equation*}
	- i e^{ -i t k^2} \int_0^t e^{i s k_1^2} 	ds  + \mathcal{O}(t^2 ( - 2 k_1 (k_2 + k_3) + 2 k_2 k_3 )).
\end{equation*}
The decorated trees for NLS are of different nature and are given for $r=2$ by
\begin{equation*}
	\begin{aligned}
	\mathcal{T}_k^2  &= \left\lbrace T_0, T_1, T_2, T_3, \, k_i \in \mathbb{Z}^d \right\rbrace, 
	\\ T_0 & =  
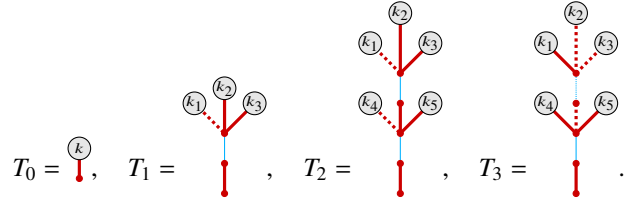
\begin{tikzpicture}[scale=0.2,baseline=-5]
		\coordinate (root) at (0,1);
		\coordinate (tri) at (0,-1);
		\draw[kernels2] (tri) -- (root);
		\node[var] (rootnode) at (root) {\tiny{$ k $}};
		\node[not] (trinode) at (tri) {};
	\end{tikzpicture} , \quad T_1 =  \begin{tikzpicture}[scale=0.2,baseline=-5]
		\coordinate (root) at (0,2);
		\coordinate (tri) at (0,0);
		\coordinate (trib) at (0,-2);
		\coordinate (t1) at (-2,4);
		\coordinate (t2) at (2,4);
		\coordinate (t3) at (0,5);
		\draw[kernels2,tinydots] (t1) -- (root);
		\draw[kernels2] (t2) -- (root);
		\draw[kernels2] (t3) -- (root);
		\draw[kernels2] (trib) -- (tri);
		\draw[symbols] (root) -- (tri);
		\node[not] (rootnode) at (root) {};
		\node[not] (trinode) at (tri) {};
		\node[var] (rootnode) at (t1) {\tiny{$ k_{\tiny{1}} $}};
		\node[var] (rootnode) at (t3) {\tiny{$ k_{\tiny{2}} $}};
		\node[var] (trinode) at (t2) {\tiny{$ k_{\tiny{3}} $}};
		\node[not] (trinode) at (trib) {};
	\end{tikzpicture}, \quad  T_2 = \begin{tikzpicture}[scale=0.2,baseline=-5]
		\coordinate (root) at (0,2);
		\coordinate (tri) at (0,0);
		\coordinate (trib) at (0,-2);
		\coordinate (t1) at (-2,4);
		\coordinate (t2) at (2,4);
		\coordinate (t3) at (0,4);
		\coordinate (t4) at (0,6);
		\coordinate (t41) at (-2,8);
		\coordinate (t42) at (2,8);
		\coordinate (t43) at (0,10);
		\draw[kernels2,tinydots] (t1) -- (root);
		\draw[kernels2] (t2) -- (root);
		\draw[kernels2] (t3) -- (root);
		\draw[symbols] (root) -- (tri);
		\draw[symbols] (t3) -- (t4);
		\draw[kernels2,tinydots] (t4) -- (t41);
		\draw[kernels2] (t4) -- (t42);
		\draw[kernels2] (t4) -- (t43);
		\draw[kernels2] (trib) -- (tri);
		\node[not] (trinode) at (trib) {};
		\node[not] (rootnode) at (root) {};
		\node[not] (rootnode) at (t4) {};
		\node[not] (rootnode) at (t3) {};
		\node[not] (trinode) at (tri) {};
		\node[var] (rootnode) at (t1) {\tiny{$ k_{\tiny{4}} $}};
		\node[var] (rootnode) at (t41) {\tiny{$ k_{\tiny{1}} $}};
		\node[var] (rootnode) at (t42) {\tiny{$ k_{\tiny{3}} $}};
		\node[var] (rootnode) at (t43) {\tiny{$ k_{\tiny{2}} $}};
		\node[var] (trinode) at (t2) {\tiny{$ k_5 $}};
	\end{tikzpicture}, \quad T_3 = \begin{tikzpicture}[scale=0.2,baseline=-5]
		\coordinate (root) at (0,2);
		\coordinate (tri) at (0,0);
		\coordinate (trib) at (0,-2);
		\coordinate (t1) at (-2,4);
		\coordinate (t2) at (2,4);
		\coordinate (t3) at (0,4);
		\coordinate (t4) at (0,6);
		\coordinate (t41) at (-2,8);
		\coordinate (t42) at (2,8);
		\coordinate (t43) at (0,10);
		\draw[kernels2] (t1) -- (root);
		\draw[kernels2] (t2) -- (root);
		\draw[kernels2,tinydots] (t3) -- (root);
		\draw[symbols] (root) -- (tri);
		\draw[symbols,tinydots] (t3) -- (t4);
		\draw[kernels2] (t4) -- (t41);
		\draw[kernels2,tinydots] (t4) -- (t42);
		\draw[kernels2,tinydots] (t4) -- (t43);
		\draw[kernels2] (trib) -- (tri);
		\node[not] (trinode) at (trib) {};
		\node[not] (rootnode) at (root) {};
		\node[not] (rootnode) at (t4) {};
		\node[not] (rootnode) at (t3) {};
		\node[not] (trinode) at (tri) {};
		\node[var] (rootnode) at (t1) {\tiny{$ k_{\tiny{4}} $}};
		\node[var] (rootnode) at (t41) {\tiny{$ k_{\tiny{1}} $}};
		\node[var] (rootnode) at (t42) {\tiny{$ k_{\tiny{3}} $}};
		\node[var] (rootnode) at (t43) {\tiny{$ k_{\tiny{2}} $}};
		\node[var] (trinode) at (t2) {\tiny{$ k_5 $}};
	\end{tikzpicture}.
\end{aligned}
\end{equation*}
where for $T_1$, $k= -k_1 + k_2 + k_3$, for $T_2$, $ k = -k_4 -k_1 + k_2 + k_3 + k_5 $ and for $T_3$,  $ k = k_4 +k_1 - k_2 - k_3 + k_5 $. An edge  $ \begin{tikzpicture}[scale=0.2,baseline=-5]
	\coordinate (root) at (0,-1);
	\coordinate (tri) at (0,0.5);
	\draw[kernels2] (root) -- (tri);
\end{tikzpicture} $ (resp. $ \begin{tikzpicture}[scale=0.2,baseline=-5]
\coordinate (root) at (0,-1);
\coordinate (tri) at (0,0.5);
\draw[kernels2,tinydots] (root) -- (tri);
\end{tikzpicture} $) corresponds to a factor $ e^{-i t k^2} $ (resp. $ e^{i t k^2} $), while an edge  $  \begin{tikzpicture}[scale=0.2,baseline=-5]
	\coordinate (root) at (0,-1);
	\coordinate (tri) at (0,0.5);
	\draw[symbols] (root) -- (tri);
\end{tikzpicture}  $ (resp. $  \begin{tikzpicture}[scale=0.2,baseline=-5]
\coordinate (root) at (0,-1);
\coordinate (tri) at (0,0.5);
\draw[symbols,tinydots] (root) -- (tri);
\end{tikzpicture}  $)  corresponds to an integral $ - i  \int^{t}_0 e^{i s k^2} \cdots d s $ (resp. $ i  \int^{t}_0 e^{-i s k^2} \cdots d s $ ). The dotted edges can be seen as taking the complex conjugate of the operator. Also, the frequencies add up to the root with a minus sign for dotted edges. Indeed, for the decorated tree $T_3$, one has
\begin{equation*}
	 k = k_4 - \ell + k_5, \quad -\ell = k_1 -k_2 -k_3
	\end{equation*}
where $ \ell $ corresponds to the node decoration of the inner nodes not connected to the root. This shows that our decorated tree formalism is very robust for encoding many iterated integrals originating from dispersive PDEs.

One can further simplify the scheme by using a nice trick. The idea is that one may know in advance the regularity a priori asked on the initial data in some Sobolev space.  If one has enough regularity on the initial data, then the computation \eqref{resonance_decomposition} can be replaced by a classical discretisation without a resonance analysis
 \begin{equation} \label{classical_decomposition}
 	\begin{aligned}
 		\int_0^t e^{i s \mathcal{L}} ds &= 		\int_0^t (1 + \mathcal{O}(s\mathcal{L})) ds = t +  \mathcal{O}(t^2\mathcal{L}).
 	\end{aligned}
 \end{equation}
Then, given a decorated tree $T$, we denote by $(\Pi^{n,r} T)(t)$ the discretisation of $ (\Pi T)(t) $ following these rules. Here, the index $n$ corresponds to the a priori regularity of the initial data and $r$  is the order of the discretisation. If we let $n$ vary, one gets an entire family of 
resonance-based schemes.
The resonance-based scheme takes the form
\begin{equation} \label{genscheme}
	 U^{n,r}_k(v,t) = \sum_{T \in \mathcal{T}_k^r} \frac{\Upsilon(T)(v)}{S(T)} (\Pi^{n,r} T)(t).
\end{equation}
 One has the following general result for the computation of the local error
\begin{theorem} \label{approxima_tree}
	For every decorated tree $ T  $ one has,
	\begin{equation*}
		\left(\Pi T - \Pi^{n,r} T \right)(t)  = \mathcal{O}\left( t^{r+1} \mathcal{L}^{r}_{{ \normalfont \text{\tiny{low}}}}(T,n) \right).
	\end{equation*}
Then, one has
	\begin{equation*}
		U_{k}^{n,r}(v,t) - U_{k}(t) = \sum_{T \in\mathcal{T}_k^r } \mathcal{O}\left(t^{r+1} \mathcal{L}^{r}_{\normalfont \text{\tiny{low}}}(T,n) \Upsilon( T)(v) \right)
	\end{equation*}
	where { there exists an operator} $\mathcal{L}^{r}_{ \normalfont \text{\tiny{low}}}(T,n)$ which can be explicitly computed embeds the necessary regularity of the solution.
\end{theorem}

{ The operator $\mathcal{L}^{r}_{ \normalfont \text{\tiny{low}}}(T,n)$ is uniquely determined by the choice performed in the discretisation meaning the splitting of the various operators $\mathcal{L}$ connected to the decorated tree $T$ into $\mathcal{L}_{\text{\tiny{low}}}$ and $\mathcal{L}_{\text{\tiny{dom}}}$. In \cite[Def. 3.11]{BS}, one has a recursive formula for computing these operators that relies crucially on the Hopf algebraic construction via a well-chosen coproduct on decorated trees.}

As an application  of the previous theorem, one can get a low regularity scheme for the KdV equation:
\begin{corollary}\label{corKdV} For the KdV equation \eqref{kdv} the general scheme~ \eqref{genscheme} takes at first order the form
	\begin{equation}
		\begin{aligned}\label{schemeKdV1}
			u^{\ell+1} &=  e^{-\tau \partial_x^3} u^\ell + \frac16 \left(e^{-\tau\partial_x^3 }\partial_x^{-1} u^\ell\right)^2 - \frac16 e^{-\tau\partial_x^3} \left(\partial_x^{-1} u^\ell\right)^2\end{aligned}
	\end{equation}
	with a local error  of order $\mathcal{O}\Big(
	\tau^2 \partial_x^2 u
	\Big)$
	and at  second-order 
	\begin{equation}
		\begin{aligned}\label{schemeKdV}
			u^{\ell+1} &=  e^{-\tau \partial_x^3} u^\ell + \frac16 \left(e^{-\tau\partial_x^3 }\partial_x^{-1} u^\ell\right)^2 - \frac16 e^{-\tau\partial_x^3} \left(\partial_x^{-1} u^\ell\right)^2\\& +\frac{\tau^2}{4} e^{- \tau \partial_x^3}\Psi\big(i \tau \partial_x^2\big)  \Big(\partial_x \Big(u^\ell \partial_x (u^\ell u^\ell)\Big)\Big)
		\end{aligned}
	\end{equation}
	with a local error  of order $\mathcal{O}\Big(
	\tau^3 \partial_x^4 u
	\Big)$ and  a suitable filter function $\Psi$ satisfying
	\begin{equation*}
		\Psi= \Psi\left(i \tau \partial_x^2 \right), \quad \Psi(0) = 1, \quad \Vert \tau   \Psi \left(i \tau \partial_x^2\right) \partial_x^2 \Vert_r \leq 1.
	\end{equation*}
\end{corollary}

With the aid of stability one can obtain a global error estimate from the local error of Theorem~\ref{approxima_tree} with the aid of Lady Windamere's fan argument \cite{H2Tri}. The necessary stability estimates in general rely on the algebraic structure of the underlying space. In the stability analysis of dispersive PDEs set in Sobolev spaces $H^r$ one classically exploits bilinear estimates of type 
	\[
	\Vert v w \Vert_r \leq c_{r,d} \Vert v \Vert_r \Vert w \Vert_r.
	\]
	The latter only holds for $r>d/2$ and thus restricts the analysis to sufficiently smooth Sobolev spaces $H^r$ with $r>d/2$. To obtain (sharp) $L^2$ global error estimates one needs to exploit discrete Strichartz estimates and discrete  Bourgain spaces in the periodic setting, see, e.g., \cite{IZ09,ORS19,ORS20} and section \ref{sectionBourgain}.

	One can study the composition and substitution of numerical methods of the form \eqref{genscheme} repeating the construction performed for B-series (see \cite{MR2657947} for a survey). The substitution is connected to the renormalisation in quantum field theory \cite{CKI}. A first work in this direction is \cite{Br23} which computes such operations in the context of regularity structures from which the schemes borrow their form.
	
Theorem~\ref{approxima_tree} has been obtained in \cite{BS} for a large class of dispersive PDEs. It computes the local error for a resonance-based scheme. Let us briefly explain how $\mathcal{L}^{r}_{\text{\tiny{low}}}(T,n)$ can be computed via Hopf algebra techniques.
One can define a coproduct on $ \mathcal{F} $ the forest formed of decorated trees in $ \mathcal{T} $. This is the free commutative monoid generated by $\mathcal{T}$. The coproduct $\Delta_{\text{\tiny{BCK}}} : \mathcal{F} \rightarrow \mathcal{F} \otimes \mathcal{F}$ is defined via admissible cuts and it is very similar to the Butcher-Connes-Kreimer coproduct \cite{Butcher72,CK} used for numerical analysis and renromalisation in quantum field theory. We will not give a precise definition of this map but just a computation on one example coming from the KdV equation. One has
\begin{equation*}
	\Delta_{\text{\tiny{BCK}}} \begin{tikzpicture}[scale=0.2,baseline=-5]
		\coordinate (root) at (0,0);
		\coordinate (tri) at (0,-2);
		\coordinate (t1) at (-1,2);
		\coordinate (t11) at (-2,4);
		\coordinate (t12) at (-3,6);
		\coordinate (t13) at (-1,6);
		\coordinate (t2) at (1,2);
		\draw[kernels2] (t11) -- (t13);
		\draw[kernels2] (t11) -- (t12);
		\draw[kernels2] (t1) -- (root);
		\draw[symbols] (t1) -- (t11);
		\draw[kernels2] (t2) -- (root);
		\draw[symbols] (root) -- (tri);
		\node[not] (rootnode) at (root) {};
		\node[not,label= {[label distance=-0.2em]below: \scriptsize  }] (trinode) at (tri) {};
		\node[not] (trinode) at (t1) {};
		\node[var] (rootnode) at (t12) {\tiny{$ k_{\tiny{1}} $}};
		\node[var] (rootnode) at (t13) {\tiny{$ k_{\tiny{2}} $}};
		\node[var] (trinode) at (t2) {\tiny{$ k_3 $}};
	\end{tikzpicture}  = \begin{tikzpicture}[scale=0.2,baseline=-5]
		\coordinate (root) at (0,0);
		\coordinate (tri) at (0,-2);
		\coordinate (t1) at (-1,2);
		\coordinate (t11) at (-2,4);
		\coordinate (t12) at (-3,6);
		\coordinate (t13) at (-1,6);
		\coordinate (t2) at (1,2);
		\draw[kernels2] (t11) -- (t13);
		\draw[kernels2] (t11) -- (t12);
		\draw[kernels2] (t1) -- (root);
		\draw[symbols] (t1) -- (t11);
		\draw[kernels2] (t2) -- (root);
		\draw[symbols] (root) -- (tri);
		\node[not] (rootnode) at (root) {};
		\node[not,label= {[label distance=-0.2em]below: \scriptsize  }] (trinode) at (tri) {};
		\node[not] (trinode) at (t1) {};
		\node[var] (rootnode) at (t12) {\tiny{$ k_{\tiny{1}} $}};
		\node[var] (rootnode) at (t13) {\tiny{$ k_{\tiny{2}} $}};
		\node[var] (trinode) at (t2) {\tiny{$ k_3 $}};
	\end{tikzpicture} \otimes \one + \one \otimes \begin{tikzpicture}[scale=0.2,baseline=-5]
		\coordinate (root) at (0,0);
		\coordinate (tri) at (0,-2);
		\coordinate (t1) at (-1,2);
		\coordinate (t11) at (-2,4);
		\coordinate (t12) at (-3,6);
		\coordinate (t13) at (-1,6);
		\coordinate (t2) at (1,2);
		\draw[kernels2] (t11) -- (t13);
		\draw[kernels2] (t11) -- (t12);
		\draw[kernels2] (t1) -- (root);
		\draw[symbols] (t1) -- (t11);
		\draw[kernels2] (t2) -- (root);
		\draw[symbols] (root) -- (tri);
		\node[not] (rootnode) at (root) {};
		\node[not,label= {[label distance=-0.2em]below: \scriptsize  }] (trinode) at (tri) {};
		\node[not] (trinode) at (t1) {};
		\node[var] (rootnode) at (t12) {\tiny{$ k_{\tiny{1}} $}};
		\node[var] (rootnode) at (t13) {\tiny{$ k_{\tiny{2}} $}};
		\node[var] (trinode) at (t2) {\tiny{$ k_3 $}};
	\end{tikzpicture} 
	+  \begin{tikzpicture}[scale=0.2,baseline=-5]
		\coordinate (root) at (0,0);
		\coordinate (tri) at (0,-2);
		\coordinate (t1) at (-1,2);
		\coordinate (t2) at (1,2);
		\draw[kernels2] (t1) -- (root);
		\draw[kernels2] (t2) -- (root);
		\draw[symbols] (root) -- (tri);
		\node[not] (rootnode) at (root) {};
		\node[not,label= {[label distance=-0.2em]below: \scriptsize  }] (trinode) at (tri) {};
		\node[var] (rootnode) at (t1) {\tiny{$ \ell $}};
		\node[var] (trinode) at (t2) {\tiny{$ k_3 $}};
	\end{tikzpicture}   \otimes \begin{tikzpicture}[scale=0.2,baseline=-5]
		\coordinate (root) at (0,0);
		\coordinate (tri) at (0,-2);
		\coordinate (t1) at (-1,2);
		\coordinate (t2) at (1,2);
		\draw[kernels2] (t1) -- (root);
		\draw[kernels2] (t2) -- (root);
		\draw[symbols] (root) -- (tri);
		\node[not] (rootnode) at (root) {};
		\node[not,label= {[label distance=-0.2em]below: \scriptsize  }] (trinode) at (tri) {};
		\node[var] (rootnode) at (t1) {\tiny{$ k_{\tiny{1}} $}};
		\node[var] (trinode) at (t2) {\tiny{$ k_2 $}};
	\end{tikzpicture}  
\end{equation*}
where $ \ell = k_1 + k_2 $. On the right part of the tensor product, we have collected the branches of the decorated tree which have been cut. On the left part, we have put the trunk. Here, $\one$ denotes the empty tree. The edges that can be cut are in blue and there are only two choices: Cutting the edge connected to the root or the one above. It is not possible to cut them simultaneously as they are on the path to the root (not admissible cut).
Why such a  cutting algorithm can be used in the context of dispersive PDEs? This comes from the computation of the dominant parts and their exact integration in the resonance analysis. Indeed, one has
\begin{equation*}
	\int_0^t e^{i s \mathcal{L}_{\text{\tiny{dom}}}} ds = \frac{e^{i t \mathcal{L}_{\text{\tiny{dom}}}} -1}{i  \mathcal{L}_{\text{\tiny{dom}}}}
\end{equation*} 
that gives two different contributions $ e^{i s \mathcal{L}_{\text{\tiny{dom}}}} $ and $1$ to plug into a new resonance problem for computing the next steps of the discretisation. Then, the contribution $1$ corresponds to a branch we detach and $ e^{i s \mathcal{L}_{\text{\tiny{dom}}}} $ to the trunk. 
One can define a new map on decorated trees $ \hat{\Pi}^{n,r}   $ which is of the form
\begin{equation*}
(\hat{\Pi}^{n,r} T)(t) = b(n,r,T,t) e^{it\mathscr{F}_{\text{\tiny{dom}}}(T)}.
	\end{equation*}
where $ \mathscr{F}_{\text{\tiny{dom}}}(T) $ is computed recursively and corresponds to the dominant resonance coming from $ \Pi^{n,r} T $. Here, $b(n,r,T,t)$ is a coefficient depending on many paremeters. In fact, one is able to write a Birkhoff type factorisation of the form
\begin{equation*}\label{Birkhoff2}
	\hat \Pi^n = \left( \Pi^n \otimes \left( \mathcal{Q} \circ \Pi^n \mathcal{A} \cdot \right)(0) \right) \hat{\Delta}_{\text{\tiny{BCK}}}.
\end{equation*}
Here, both $ \Pi^n $ and $ \hat{\Pi}^{n} $ are defined on a larger set of decorated trees where the order of the discretisation is now part of the root decoration. The coprodcut $\hat{\Delta}_{\text{\tiny{BCK}}}$ is an extension of $ \Delta_{\text{\tiny{BCK}}} $ closer in spirit to the recentering coproduct introduced in  \cite{reg,BHZ}. The map $ \mathcal{A} $ is a recursive map called antipode that performs the various non-admissible cuts. The map $\mathcal{Q}$ is a projector that keeps only the mode $1$ removing all the terms of the form $e^{i P(k_1,...,k_n)}$. This Hopf algebra approach has been pushed forward in \cite{ABBMS23} by proposing a forest formula for the resonance-based schemes to find symmetric schemes. It is an open problem to get very explicit formulae for these schemes although we know that they depend on various lower and dominant parts that can be computed explicitly.

Before moving to structure-preserving schemes,  one can generalise the resonance approach to non-polynomial nonlinearities of type $$f(u) g(\overline u)$$ for smooth functions $f$ and $g$. { We allow $f$ and $g$ to be equal to one. For example, in KdV, one does not have a dependency in $\bar{u}$.} Indeed, if one gets   \[
 e^{i s\mathcal{L}} v + A(v,s) 
\]
where $ A(v,s) $   a linear combination of iterated integrals. Then, we  perform a Taylor expansion around the point $ e^{is \mathcal{L}}v $:
\begin{equation*}
	f(e^{i s\mathcal{L}} v + A(v,s)) = \sum_{m \leq r} \frac{A(v,s)^m}{m!} f^{(m)}(e^{i s\mathcal{L}} v ) + \mathcal{O}(A(v,s)^{r+1}).
\end{equation*}
Performing the same type of computations on $g$, one obtains
$$
\frac{A(v,s)^m}{m!} f^{(m)}(e^{i s\mathcal{L}} v ) 
\frac{\overline{A(v,s)}^n}{n!} g^{(n)}({e^{- i s\mathcal{L}}\overline v}) .
$$
We need to pull the oscillatory phases $ e^{\pm is\mathcal{L}} $ out of $f$ and $g$. This is achieved via expansions of the form
\begin{equation*}
	f(e^{is \mathcal{L}}v) = \sum_{\ell \leq r} \frac{s^{\ell}}{\ell!} e^{is \mathcal{L}} \mathcal{C}^{\ell}[f,\mathcal{L}](v) + \mathcal{O}(s^{r+1} \mathcal{C}^{r+1}[f,\mathcal{L}](v))
\end{equation*}
where $\mathcal{C}^{\ell}[f,\mathcal{L}]$ denote nested commutators which in general require less regularity than powers of the full operator $\mathcal{L}^\ell$. Then, one can perform the resonance analysis in Fourier space or continue directly in physical space by using again the same commutators. 
This approach is fairly general and allows us to cover the time dynamics of a large class of equations, including parabolic
and hyperbolic problems, as well as dispersive equations, up to arbitrarily high
order on general domains (see \cite{RS,ABBS22}).

With the naive choices performed for deriving our schemes, it is quite likely that the scheme does not preserve any structure. There are various places where one can improve the discretisation:
\begin{itemize}
	\item One can perform a different Duhamel iteration. Indeed, one has the general iteration for $ s \in [0,t] $
	\begin{equation*} 
		\begin{aligned}
			I(k,u,s,t) 
			 & = e^{i (t-s) k^3 }   u_k(s)  + i k \sum_{\substack{k_1,k_2 \in \mathbb{Z}\\k_1+k_2 = k} } i e^{ i t k^3}  \int_s^t e^{-i\tilde{s} k^3} 	u_{k_1}(\tilde{s}) 	u_{k_2}(\tilde{s}) d\tilde{s}, 
		\end{aligned}
	\end{equation*}  
We can take a weighted sum of these Duhamel iterations
\begin{equation*} \label{mid_point}
	u_k(t) = \frac{1}{2} \left(  I(k,u,0,t) + I(k,u,s,t) \right).
\end{equation*}
which produced the midpoint iteration. One has to work with an implicit scheme as the value of both $u_k(s)$ and $u_k(0)$ are needed. Let us mention that this trick could  be used outside numerical analysis as it gives a certain degree of freedom in the way one iterates Duhamel's formula.
\item Another degree of freedom is when we perform  Taylor expansions. Instead of doing them around $0$, we can take an interpolation. We consider $ r +1 $ distinct interpolation points $0\leq a_0<a_1<\dots<a_r\leq 1$ which are symmetrically distributed such that $a_{j}=1-a_{r-j},j=0,\dots, r$. Let us denote the corresponding nodal polynomials by $ p_{j,r} $ such that
\begin{equation*}
	p_{j,r}(a_m t) = \delta_{j,m}.
\end{equation*}
Then, we define the following approximation
\begin{equation*}
	\tilde{p}_r( f,s) = \sum_{j=0}^r f(a_j)  p_{j,r}(s), \quad f(a_j) =  e^{i a_j t \mathcal{L}_{\text{\tiny{low}}}}.
\end{equation*}
\end{itemize}

With the previous ideas coming from \cite{MS22,ABBMS23}, one can derive a low regularity symmetric scheme for the KdV equation. We recall briefly what is the definition of a symmetric scheme. It is defined by considering its adjoint method: For a given method $u^{\ell}\mapsto u^{\ell+1}=\Phi_{t}(u^{\ell})$ its adjoint method is defined as $\widehat{\Phi}_{t}:=\Phi^{-1}_{-t}$.
\begin{definition} The method $\Phi_{t}$ is called \textit{symmetric} if $\Phi_{t}=\widehat{\Phi}_{t}$.
\end{definition}

\begin{corollary}\label{corKdV} For the KdV equation \eqref{kdv} the general midpoint scheme~ \eqref{genscheme} takes at first order the form
	\begin{equation}
		\begin{aligned}\label{schemeKdV1}
			u^{\ell+1} &=  e^{-\tau \partial_x^3} u^\ell + \frac{1}{24} \left(e^{-\tau\partial_x^3 }\partial_x^{-1}  u^{\ell} + \partial_x^{-1} u^{\ell+1} \right)^2 \\ & - \frac{1}{24} e^{-\tau\partial_x^3} \left(\partial_x^{-1}  u^{\ell} + e^{\tau \partial_x^3} \partial_x^{-1} u^{\ell+1} \right)^2\end{aligned}
	\end{equation}
	with a local error  of order $\mathcal{O}\Big(
	\tau^2 \partial_x^2 u
	\Big)$
	at  first-order
	and
	with a local error  of order $\mathcal{O}\Big(
	\tau^3 \partial_x^4 u
	\Big)$ at second order. 
\end{corollary}

Symmetric schemes have been understood in a great degree of generality in \cite{ABBMS23}.   One may  also want to get low regularity {\em symplectic} schemes for Hamiltonian systems. The Poisson structure of the KdV equation in Fourier variables for instance takes the form (see, e.g.,  \cite{BernierGrebert})
$$
\{P,Q \}(u) =\sum_{k\in \mathbb{Z}^\ast}(\partial_{u_{-k}}P)(2i\pi k) (\partial_{u_k}Q)
$$
and it is an interesting question in how far we can reproduce this structure at the discrete (numerical) level.

For doing so, one may want to exploit another degree of freedom which is the splitting of the differential operator $ \mathcal{L} $ into
 \begin{equation*}
 	\mathcal{L} = \mathcal{S}_{\text{\tiny{dom}}} + \mathcal{S}_{\text{\tiny{low}}}
 	\end{equation*}
 where the splitting has to preserve the symplectic structure of the equation. In the case of the KdV equation one does not have the choice for the first iterated integral as one performs directly an exact integration that preserves the symplectic structure.
For other equations such as for instance cubic NLS,  one changes the splitting:
 \begin{equation*}
 	\mathcal{S}_{\text{\tiny{dom}}}(k,k_1)  = - 2 k_1 k, \quad 	\mathcal{S}_{\text{\tiny{low}}}(k_2,k_3)  =  2 k_2 k_3.
 \end{equation*}
What is crucial is that now, one has 
\begin{equation} \label{symmetry_symplectic}
	e^{i s 	\mathcal{S}_{\text{\tiny{dom}}}(k,k_1) } e^{i s 	\mathcal{S}_{\text{\tiny{low}}}(k_2,k_3) } = \overline{ e^{i s 	\mathcal{S}_{\text{\tiny{dom}}}(k_2,k_3) } e^{i s 	\mathcal{S}_{\text{\tiny{low}}}(k_1,k_2) }}
\end{equation}
 which corresponds to the symmetry coming from the symplectic nature of the NLS equation. Then in the descritisation, one has to perform almost identical steps on $	e^{i s 	\mathcal{S}_{\text{\tiny{dom}}}(k,k_1) }$ and $  e^{i s 	\mathcal{S}_{\text{\tiny{low}}}(k_2,k_3) } $ to respect \eqref{symmetry_symplectic}. The discretisation is limited to the dimension one ($d=1$) (see \cite{MS22} and \cite{AGB24} in the SPDE context). It is an open problem to get low regularity symplectic  schemes for any dimension and a large class of dispersive equations. The main restriction is on the dimension as one cannot go back to physical  space for $ \frac{1}{kk_1} $ with usual operators when $d$ is larger than one.

Let us finish this section  by providing some perspectives. As we have already mentioned, structure-preserving low regularity   schemes are quite challenging to get. However, many recent techniques have allowed some significant progress in this direction.
On the combinatorial side, the structures proposed for describing the schemes have been used in a more analytical context. 

Indeed, Poincaré-Dulac normal forms  for dispersive PDEs have been derived  in \cite{Br24} with the same decorated trees and an arborification map. This gives a new derivation of this normal form \cite{oh1} often used in the context of dispersive PDEs.
This derivation is also very natural as it is reminiscence of structures that have been used in a more general context such as mould calculus introduced by Jean Ecalle \cite{Ecalle1,Ecalle2} for dynamical systems. The same arborification could be found there \cite{EV04,FM} but also in teh context of rough analysis in \cite{U10}. This formalism can even explain some of the most recent progress for dispersive PDEs with random initial data. 
In \cite{BT24}, one can compute hidden cancellations with a similar arborification for dispersive PDEs. It covers cancellations coming 
 from Wave turbulence \cite{DH23,DH2301} and the proof of the invariance of the
Gibbs measure under the dynamics of the three-dimensional cubic wave equation \cite{BDNY24}. Let us mention low regularity schemes in the context of Wave turbulence have been derived in \cite{ABBS24}.
The decorated formalism could also cover the recent theory of random tensors \cite{DNY22} that needs to be compared with some word formalism developped in \cite{WordSeries}.

The combinatorial formalism proposed for low regularity schemes for dispersive PDEs is very robust and it connects to various other fields
such as dynamical systems, noncommutative geometry, numerical analysis with B-series, renormalisation in quantum field theory, rough analysis, singular SPDEs.
 Then, some ideas can be borrowed from these fields for organising  computations connected to dispersive PDEs but also beyond.

\section{Error analysis at low regularity with discrete Bourgain spaces}
\label{sectionBourgain}
In this section, we shall survey the main ideas which allows to perform error estimates at low regularity by taking the KdV equation
\eqref{kdv} as the main example. We shall focus on error estimates in $L^2$ and thus restrict ourselves to solutions with positive regularity.

A flexible tool for the study of the local well-posedness of dispersive partial differential equations at low regularity is the use  of Bourgain spaces
introduced in \cite{Bour93b}.
We refer for more details to the books \cite{Tao06}, \cite{Linares}.

{
To make use of these spaces at the discrete -- numerical -- level we make the following  observation.
\begin{remark}\label{rem:Pi}
In order to exploit discrete Bourgain spaces in numerical analysis, we have to introduce frequency cut-off in our numerical schemes.   This is crucial in our estimates to be able to push down the regularity assumptions on our solution. Indeed in the discrete case, the dispersion relation
  is periodic  of period $2 \pi/ \tau$ with respect to $\sigma$, this create additional frequency interactions,  compared to the continuous case,  which could not
  be controlled without loosing  derivatives. The presence of the cut-off allows us to have at the discrete level exactly the 
  same type of frequency interactions as in the continuous case and in particular to be able to use the factorization property \eqref{factorization}
  which is crucial to recover the derivative in KdV.
  
  The introduction of a frequency cut-off $ \Pi_{\tau}$ in the scheme (see \eqref{bourgNL} for its precise structure) does not change the resonance analysis given in the previous section. Indeed, the filtered scheme can be seen as a resonance based discretisation of the corresponding filtered PDE, i.e., 
  $$
  \partial_{t} u^\tau + \partial_{x}^3 u^\tau + \Pi_{\tau}(  \Pi_{\tau } u^\tau \Pi_{\tau} \partial_{x}  u)=0
  $$
  in case of KdV.
  
  It is also interesting to mention that   frequency cut-offs  are natural in the numerical discretisation of PDEs (thinking for instance about spectral methods) as we can not reproduce all frequencies on the computer. Our approach is to  choose this frequency cut-off in such a way that we can exploit discrete Bourgain spaces and thus optimise the error estimates in sense of regularity.
\end{remark}
}
 A tempered distribution $u(t,x)$ on $\mathbb{R} \times \mathbb{T}$ belongs to the Bourgain space $X^{s, b}$ adapted to the Airy equation
(the linear part of \eqref{kdv}) if its following norm is finite
\begin{equation*}
\|u\|_{X^{s, b}}=\left(\int_{\mathbb{R}}\sum_{k \in \mathbb{Z}}\left(1+ |k|\right)^{2s}\left(1+| \sigma -  k^3|\right)^{2b}|\widetilde{u}\left(\sigma, k\right)|^2 \,d\sigma\right)^{\frac{1}{2}},
\end{equation*}
where $\widetilde{u}$ is the space-time Fourier transform of $u$:
$$
\widetilde{u}(\sigma, k)= \int_{\mathbb{R}\times \mathbb{T}} e^{-i \sigma t - i k x} u(t,x) \, d t  d x.
$$

This norm can be also written as
$$ \|u \|_{X^{s,b}}= \| e^{t \partial_{x}^3} u \|_{H^b_{t} H^s_{x}}$$
where $ H^{b}_{t} H^s_{x}$ is  the  space time Sobolev space, whose norm is defined by 
$ \|f\|_{H^b_{t} H^s_{x}} := \| (1+| \sigma|)^b (1 +|k|)^s \tilde f\|_{L^2(\mathbb{R} \times \mathbb{T})}.$

Elementary properties of the Bourgain spaces are the following. 

 For $\eta\in \mathcal{C}^\infty_{c}(\mathbb{R})$,   we have that
 \begin{eqnarray}
 \label{bourg1}
&  \| \eta (t) e^{- t \partial_{x}^3 } f \|_{X^{s, b}} \lesssim_{\eta} \|f\|_{H^s}, \quad s \in \mathbb{R}, \, f \in H^s(\mathbb{T}), \\
  \label{bourg2}
   \label{bourg3}
&  \| \eta({t \over T}) u \|_{X^{s,b'}} \lesssim_{\eta, b, b'} T^{b-b'} \|u\|_{X^{s,b}}, \quad s \in \mathbb{R},  -{1 \over 2} <b' \leq b <{ 1 \over 2}, \, 0< T \leq 1, \\
\label{bourg4}
& \left\|  \eta (t) \int_{0}^t e^{- (t-s) \partial_{x}^3} F(s) \, ds \right\|_{X^{s, b}} \lesssim_{\eta, b} \|F\|_{X^{s, b-1}}, \quad s \in \mathbb{R}.
 \end{eqnarray}

 In order to solve by fixed point an  equation under the form 
 $$ \partial_{t} u  + \partial_{x}^3 u = B[u, u]$$
 where $B[\cdot, \cdot ]$ is a bilinear term, one can notice that a solution in short times solves the modified Duhamel Formula
 $$ u(t)= \eta(t) e^{-t \partial_{x}^3} u_{0} + \eta(t) \int_{0}^t e^{-(t-s) \partial_{x}^3}  \eta( {s \over \delta})B[  u,  u ] (s) \,ds.$$
 
 By using the above elementary properties, local well-posedness for data in $H^s$  would follow from the Banach fixed point if  the bilinear estimate
 \begin{equation}
 \label{bilinear} \| B[u, v]  \|_{X^{s, b'} } \lesssim \| u \|_{X^{s,b}} \|v\|_{X^{s,b}}
 \end{equation}
 for any $0>b'>b-1$, $b >1/2$ holds.
 
 Unfortunately, for the periodic KdV equation, the bilinear term being $B[u,u]= - \partial_{x}( u^2/2)$, in order to recover one full derivative, one
 needs to take $b=1/2$, $b'=-1/2$ \cite{KPV,CKSTT} which is a borderline case for which \eqref{bourg4} for example does not hold.
 To deal with  a space which keeps  the good property of being a subspace
 of $ \mathcal{C}(\mathbb{R}, H^s)$, we can use  the framework of \cite{CKSTT} and  the smaller space
 $X^s$, which has the same scaling properties in time as $X^{s, {1 \over 2} }$,  defined by the following norm:
 \begin{equation}
 \label{Xsdef}
  \|u\|_{X^s} = \|u \|_{X^{s, {1 \over 2 }}} + \|\langle k \rangle^s \tilde u \|_{l^2(k)L^1(\sigma)}.
  \end{equation}
  We  define more precisely   $X^s$ as the  space of  space-time tempered distributions  such that $ \tilde{u}(\sigma, 0)= 0$ and
   the above norm is finite.
   The norm  $\|\langle k \rangle^s \tilde u \|_{l^2(k)L^1(\sigma)}$ is defined by 
   $$ \|\langle k \rangle^s \tilde u \|_{l^2(k)L^1(\sigma)}
   =\left( \sum_{k} \left|\int_{\mathbb{R}}| \tilde u (k , \sigma) |\, d \sigma \right|^2\right)^{1 \over 2}.$$
   
  A localized version $X^s(I)$ for $I$ an interval of time is obtained  by setting 
  $$ \|u\|_{X^{s}(I)} = \inf\{\|\overline{u} \|_{X^{s}}, \, \overline{u}|_{I} = u  \}.$$
  A well-posedness result for \eqref{kdv} at the $L^2$ level reads:
  \begin{theorem}
For every $T>0$ and $u_{0}\in L^2$, $\int_{\mathbb{T}}u_{0}= 0$,  there exists a unique solution $u$ of \eqref{kdv} such that $u \in  X^0(0,T)$. Moreover, if $u_{0} \in H^{s_{0}}$, $s_{0}> 0$, then $ u\in X^{s_{0}}(0,T)$.
\end{theorem}
Note that we  have $X^0(0,T) \subset \mathcal{C}([0, T], L^2)$.

We can now introduce the discrete (in time) counter part of the above spaces.

 For sequences of functions $(u^{n}(x))_{n \in \mathbb{Z}},$ we define the Fourier transform $\widetilde{u^{n}}(\sigma, k)$ by
$$
\mathcal F_{n,x}(u^n)(\sigma,k) =\widetilde{u^{n}} (\sigma, k)= \tau \sum_{m \in \mathbb{Z}} \widehat{u^{m}}(k) \, e^{i m \tau \sigma}, \quad \widehat{u^{m}}(k)= {1 \over 2\pi} \int_{-\pi}^\pi u^{m}(x) \,e^{-i k x}\,d x.
$$
Parseval's identity then reads
\begin{equation}\label{parseval}
\| \widetilde{u^{n}}\|_{L^2l^2}= \|u^{n}\|_{l^2_{\tau}L^2},
\end{equation}
where
$$
\| \widetilde{u^{n}}\|_{L^2l^2}^2 = \int_{-{\pi \over \tau}}^{\pi\over \tau} \sum_{k \in \mathbb{Z}}
|\widetilde{u^{n}}(\sigma, k)|^2 \,d \sigma, \quad
\|u^{n}\|_{l^2_{\tau}L^2}^2 = \tau \sum_{m \in \mathbb{Z}} \int_{-\pi}^\pi  |u^{m}(x)|^2 \,d x.
$$
We  define the discrete Bourgain spaces $X^{s,b}_\tau$ for $s\ge 0$, $b\in\mathbb R$, $\tau>0$ by
\begin{equation}\label{norm2}
\| u^n \|_{X^{s,b}_{\tau}} = \left\| \langle k \rangle^s \langle  d_{\tau}(\sigma + k^3)  \rangle^b \widetilde{u^n}(\sigma, k)  \right\|_{L^2l^2},
\end{equation}
where  $d_{\tau}(\sigma)=\frac{e^{i \tau \sigma} - 1}\tau$.
Note that $d_{\tau}$ is $2\pi/\tau$ periodic and that uniformly in $\tau$, we have $|d_{\tau}(\sigma)| \sim | \sigma |$ for $|\tau \sigma | \leq \pi$. Since $|d_{\tau}(\sigma)| \lesssim \tau^{-1}$, we also have
%
  that the discrete spaces satisfy the embeddings
\begin{equation}\label{embdisc1}
\|u^{n}\|_{X^{0, b}_{\tau}} \lesssim { 1 \over  \tau^{b-b'}} \|u^{n}\|_{X^{0, b'}_\tau}, \quad b \geq b'
\end{equation}
note that this estimate is not uniform in the time discretization parameter.

With this definition, we have the discrete counter part of the properties  \eqref{bourg1}, \eqref{bourg2}, \eqref{bourg3}:
\begin{lemma}\label{bourgainfaciled}
For $\eta \in \mathcal{C}^\infty_{c}(\mathbb{R})$ and $\tau\in(0,1]$, we have that
\begin{align}
\label{bourg1d} &\| \eta(n \tau)  e^{- n \tau \partial_{x}^3} f\|_{X^{s,b}_{\tau}} \lesssim_{\eta, b} \|f\|_{H^s}, \quad s \in \mathbb{R}, \, b \in \mathbb{R}, \, f \in H^s, \\
\label{bourg2d} &\| \eta(n \tau)  u^{n}\|_{X^{s,b}_{\tau}} \lesssim_{\eta, b} \|u^{n}\|_{X^{s,b}_{\tau}}, \quad s \in \mathbb{R}, \, b \in \mathbb{R} , \, u^{n} \in X^{s,b}_{\tau},\\
\label{bourg3d} &\left\| \eta\left(\frac{n\tau}T \right) u^{n} \right\|_{X^{s,b'}_{\tau}} \lesssim_{\eta, b, b'} T^{b-b'} \|u^{n}\|_{X^{s,b}_{\tau}}, \quad s \in \mathbb{R},  -{1 \over 2} <b' \leq b <{ 1 \over 2},\, 0< T = N \tau  \leq 1, \, N \geq 1.
\end{align}
In addition, for
$$
U^{n}(x)= \eta(n \tau) \tau \sum_{m=0}^n  e^{- ( n-m ) \tau \partial_{x}^3}  u^{m}(x),
$$
we have
\begin{equation}
\label{bourg4d}\|U^{n}\|_{X^{s,b}_{\tau}} \lesssim_{\eta, b} \|u^{n}\|_{X^{s, b-1}_{\tau}}, \quad s \in \mathbb{R}, \, b>1/2.
\end{equation}
\end{lemma}
We stress that all the above estimates are uniform in $\tau$.

The next step that we shall need in order to handle the KdV equation is to adapt  \eqref{bourg4} in the case 
$b=1/2$. We thus  define the discrete counterpart $X^s_{\tau}$ of the $X^s$ space.
We say that a sequence of function $(u^{n}(x))_{n} \in l^2_{\tau}L^2$ such that 
$ \int_{\mathbb{T}} u^n=0, \, \forall n$ is in $X^s_{\tau}$ for $s \geq 0$ if 
the following norm is finite
$$ \|u^{n}\|_{X^s_{\tau}}= \|u^{n}\|_{X^{s, {1 \over 2}}_{\tau}} + \|\langle k \rangle^s \tilde u( \sigma, k) \|_{l^2(k)L^1(\sigma)}$$
and in the same way,  we also define $Y_{\tau}^s$ by 
 $$ \|F^{n}\|_{Y^s_{\tau}}= \|F^{n}\|_{X^{s, -{1 \over 2}}_{\tau}}  +  \left\|{\langle k \rangle^s \over \langle d_{\tau}( \sigma + k^3\rangle)} \widetilde{F^{ n}}(\sigma, k) \right\|_{l^2(k)L^1(\sigma)}.$$

 \begin{lemma}
 \label{lemb1/2}
  We have the following properties:
 \begin{enumerate}
 \item We have the embedding $X^s_{\tau} \subset l^\infty(\mathbb{Z}, H^s(\mathbb{T}))$:
 \begin{equation}
 \label{sob1/2}
  \sup_{n} \| u^{n}\|_{H^s(\mathbb{T})} \lesssim \|u^{n} \|_{X^s_{\tau}}, \quad s \in \mathbb{R}, \, (u^{n})_{n} \in X^s_{\tau};
  \end{equation}
 \item  Let us define for $(u^{n})_{n} \in Y^s_{\tau}$,  and $\eta \in \mathcal{C}^\infty_{c}(\mathbb{R})$
 \begin{equation}
 \label{duhabourg}
U^{n}(x): = \eta(n \tau) \tau \sum_{m=0}^n  e^{- ( n-m ) \tau \partial_{x}^3}  u^{m}(x),
\end{equation}
then, we have
\begin{equation}
\label{bourg1/2}\|U^{n}\|_{X^{s}_{\tau}} \lesssim_{\eta} \|u^{n}\|_{Y^{s}_{\tau}}, \quad s \in \mathbb{R}.
\end{equation}
 \end{enumerate}
 The above estimates are uniform for $\tau \in (0, 1]$.
 \end{lemma}
 
 The last step in order to use these spaces for the analysis of numerical scheme at low regularity is to establish
 bilinear estimates like \eqref{bilinear}. A convenient one for the KdV equation is the following:
 
 For every $s \geq 0$, there exists $C>0$ such that for every $(u^n)_{n}$, $(v^n)_{n} \in X^{s}_{\tau}$, we have the estimate
 \begin{equation}
 \label{bourgNL} \| \partial_{x} \Pi_{\tau} \left( \Pi_{\tau} u^n\,  \Pi_{\tau}v^n \right) \|_{Y^{s}_{\tau}}
  \leq C \left(\|u^n \|_{X^{s, {1 \over 2} }_{\tau}} \|v^n \|_{X_{\tau}^{s, {1 \over 3}}} + \|v^n\|_{X^{s, {1 \over 2} }_{\tau}} \|u^n \|_{X_{\tau}^{s, {1 \over 3}}}
  \right)
  \end{equation}
  
  where $\Pi_{\tau}$ is a projection on frequencies $|k|$  smaller than $\tau^{- {1 \over 3}}$, see also Remark \ref{rem:Pi}

  Once these estimates are established, the convergence analysis at the $L^2$ level  of a numerical scheme for KdV, for example
  a fitered version of  \eqref{schemeKdV1} follows the following steps:
  \begin{itemize}
  \item Since the filtered scheme can be seen as a scheme for a filtered PDE, 
\begin{equation}\label{Fpde}
 \partial_{t} u^\tau + \partial_{x}^3 u^\tau + \Pi_{\tau}(  \Pi_{\tau } u^\tau \Pi_{\tau} \partial_{x}  u)=0,
 \end{equation}
  the first step is to establish an error estimate for 
  $$ \|u  - u^{\tau}\|_{L^2}.$$
  This can be performed at the continuous level by using the continuous Bourgain spaces.
  \item  The second step is to estimate the  consistency error of the scheme
 $\mathcal{E}_{n}$ in the discrete Bourgain space $Y^0_{\tau}$.
 \item The last step is to estimate $e^{n}= u^n- u^\tau(t_{n})$  in the Bourgain space $X^0_{\tau}$ which allows us to get
 a convergence result  for $  \sup_{n }\|e^{n}\|_{L^2}.$
    This is done by writing the discrete Duhamel formula and by using  \eqref{bourg1/2} and \eqref{bourgNL}.
  \end{itemize}
  
  We refer for example to \cite{ORS20},  \cite{ORS21}, \cite{ROS24} \cite{RS22}, \cite{RS22} for the implementation of this strategy
  in the analysis of various time discretizations of dispersive partial differential equations.

\end{document}